\documentclass[11pt]{article}

\newcommand{\bemm}{\begin{multline}}
\newcommand{\enm}{\end{multline}}
\newcommand{\beq}{\begin{equation}}
\newcommand{\eeq}{\end{equation}}
\newcommand{\beqa}{\begin{eqnarray}}
\newcommand{\eeqa}{\end{eqnarray}}
\newcommand{\ba}{\begin{array}}
\newcommand{\ea}{\end{array}}
\newtheorem{teo}{Theorem}

\newcommand{\be}{\begin{equation}}

\newcommand{\ee}{\end{equation}}
\newcommand{\bt}{\begin{teo}}
\newcommand{\et}{\end{teo}}

\newcommand{\xn}{x_1,\dots,x_n}
\newcommand{\C}{\mathbb{C}}

\newcommand{\ic}{\mathcal{I}}
\newcommand{\is}{I_{S}}

\newcommand{\al}{{\gamma}}

\newcommand{\R}{\mathbb{R}}

\newcommand{\vv}{{\bf V}}

\usepackage{amsmath}
\usepackage{amssymb}
\begin{document}

\title{{Time-Reversibility in 2-Dim Systems  }\\
{Valery G. Romanovski}\\
{Center for Applied Mathematics and
Theoretical Physics,\\
 University of Maribor,
 Krekova 2, SI-2000 Maribor,\\
 Slovenia }\\
{e-mail: valery.romanovsky@uni-mb.si }
\date{}
}



\maketitle


\begin{abstract}
We present an algorithm for finding 
 all time-reversible systems within a
given family of  2-dim  systems of ODE's whose right-hand sides are
polynomials. We also study an interconnection of time-reversibility 
and invariants of a subgroup  of $SL(2,\C)$. 
\end{abstract}

\section{ Introduction}
Time-reversible symmetry arises frequently in many studies in
classical and quantum mechanics. A 
general introductory account and a survey of state of  the art with
regards to time-reversible symmetry in dynamical systems theory and its applications to some physical problems can be found in  \cite{Lamb}.
Another recent stream of research is devoted  to investigation 
of isochronicity and linearizability of time-reversible polynomial systems
(see, e.g., 
\cite{Cairo,Chavarriga-1,Chavarriga-2,MMV,Valery-1} and  references therein). 

In the context  of dynamical systems described by ODE's of the form
\begin{equation} \label{fn}
\frac {d{\bf z}}{dt} = F({\bf z})\quad ({\bf z}\in \Omega),
\end{equation}
where $F:\Omega\mapsto T\Omega$ is a vector field and $\Omega$ is a manifold, 
by
time-reversible symmetry  we mean an invertible map
$R:\Omega\mapsto\Omega$, such that
 \be \label{xyn}
 \frac {d(R{\bf
z})}{dt}  =-F(R{\bf z}). \ee

In this paper we  investigate probably the  simplest
 time-reversible symmetries,
namely, we  limit our study to the case of 2-dim systems, that
is, ${\bf z} =(x,y) \in \mathbb{C}^2$ (or $\R^2$),  and to the
symmetries of the form
\begin{equation} \label{xy}
R : x \mapsto \al y, \ y\mapsto \al^{-1} x
\end{equation}
with $\al \in \mathbb{C}$ (or $\R$) (note, that $R$ defined by (\ref{xy}) is an involution).
Using methods of computational algebra we   show an  interconnection of the physically important phenomena of time-reversibility, involution and group invariants in dynamical systems described by autonomous polynomial two-dimensional systems of ODE's.

\section{Preliminaries}

We say that a straight line $L$ is an \emph{axis of symmetry}   
of a real autonomous two-dimensional  system of ODE's  
 if as
point--sets (ignoring the sense of the parametrization by 
time $t$) the orbits
of the system are symmetric with respect to the line $L$.
 There are two types of symmetry of a real system
with respect to a line $L$: \emph{mirror symmetry}, 
meaning that when the phase portrait is reflected in the line $L$
it is unchanged; and \emph{time--reversible symmetry}, meaning that 
when the phase portrait is reflected in the line $L$ and then the sense of
every trajectory is reversed (corresponding to a reversal of time), the
original phase portrait is obtained.
The symmetry (\ref{xy}) can be considered as a generalization of the second 
type of   reflection  to the case of complex two-dimensional systems
  in the following sense.
 
Consider  a parametric family of real systems
 \be \label{s.uv}
\dot u=U(u,v), \quad \dot v=V(u,v).
 \ee
Introducing  the complex structure on the plane $(u,v)$ by setting
$x=u+iv$ we obtain from (\ref{s.uv})   the equation
 \be \label{s.P} \dot x=P(x,\bar x), \qquad (P=U+i V).
 \ee
 We  add to this equation its complex conjugate to obtain the
 system
\[\dot x=P(x,\bar x), \ \dot {\bar x}= \overline{P(x,\bar x)}. \] 
Let
us now consider $\bar x$ as a new variable $y$ and allow the
 parameters of the second equation to be arbitrary.
 Then  (\ref{s.P}) yields  the complex system
\begin{equation} \label{eqgen}
\dot x = P(x,y), \ \ \dot y = Q(x,y).
\end{equation}

Let $a$ denote the
vector of coefficients of the polynomial $P(x,\bar x)$ in
(\ref{s.P}), arising from the real system (\ref{s.uv}) by setting $x = u + iv$.
It is easy to see that if
 $a = \pm \bar a$ (meaning that either all the coefficients are real or all 
are pure imaginary),
then the $u$--axis is an axis of symmetry of the real
system (\ref{s.uv}) and of the corresponding complex differential equation
(\ref{s.P}).
Thus,  $u$--axis is an axis of symmetry for (\ref{s.P}) if 
\begin{equation}
\label{crev1} 
P(\bar x,x) = -\overline{P(x,\bar x)} 
\end{equation} 
(the case $a = - \bar a$), or if 
\begin{equation} \label{crev2}
P(\bar x,x) = \overline{P(x,\bar x)} 
\end{equation}
(the case $a = \bar a$). We now observe that if condition (\ref{crev1}) is
satisfied then under the change 
$ x \to \bar x, \  \bar x \to
x
 $
{equation} 
(\ref{s.P}) is transformed to its negative, 
\begin{equation} \label{cneg} 
\dot
x = -P(x,\bar x), 
\end{equation} 
and if condition (\ref{crev2}) holds then (\ref{s.P}) is unchanged. Thus
condition
(\ref{crev2}) means that the system is reversible with respect to 
reflection
across the $u$--axis (i.e., the transformation does not change the 
system)
while condition (\ref{crev1}) corresponds to time--reversibility
 with 
respect to the same transformation.

If the line of reflection is not the $u$--axis but a distinct line $L$ then we
can apply the rotation $x_1=e^{-i\varphi}x$ through an appropriate angle $\varphi$ to
make $L$ the $u$--axis. In the new coordinates we have
\[
\dot x_1 = e^{-i\varphi} P(e^{i\varphi} x_1, e^{-i\varphi} \bar x_1)\,.
\]
By the discussion in the paragraph following (\ref{cneg}) this system is 
time--reversible with respect to
the line ${\rm Im}\, x_1=0$ if (\ref{crev1})
holds, meaning that
\[
e^{i\varphi} \overline{ P(e^{i\varphi} x_1, e^{-i\varphi}\bar x_1)}
               =
 - e^{-i\varphi}  P(e^{i\varphi}\bar x_1, e^{-i\varphi} x_1).
\]
Hence, reverting to the variable $x$, (\ref{s.P}) is time--reversible when 
there exists a $\varphi$ such that 
\begin{equation}
\label{crr} 
e^{2i\varphi} \overline{P(x,\bar x)} = - P(e^{2i\varphi}\bar x,
e^{-2i\varphi}x). 
\end{equation}
We see that in the particular case when   in (\ref{xyn}) 
$ F=(P,Q),$ 
$\al = e^{2i\varphi}$, $y = \bar x$, and
$Q = \bar P$  the equality (\ref{xyn}) is equivalent to (\ref{crr}).
Thus, the second kind  of reflection with respect to  a line
described at the beginning of the section  can be considered as a 
particular case of symmetry (\ref{xy}). 

In fact we have shown   that  if system (\ref{eqgen}) is the
complexification of a real system (\ref{s.uv}) and it admits a
symmetry (\ref{xy}) with
 $\al_0=e^{2i\phi_0}$ then the line $v =u\, \tan \phi_0$ is
the line of symmetry of the trajectories (as point-sets) 
 of the real  system (see
 \cite{Sib1} for more  details).


\section{ The set of time-reversible systems}
 Direct calculation shows that the
system (\ref{eqgen}) is time--reversible with respect to a
transformation (\ref{xy}) if and only if for some $\al$
\begin{equation} \label{crr1}
\al Q(\al y, x/\al) = -P(x, y), \quad  \al Q(x, y) = -P(\al y,
x/\al)\,.
\end{equation}
We will limit our study to the case when $P(x,y)$, $Q(x,y)$ in
(\ref{eqgen}) are polynomials.
 Without loss of generality we can write such   systems in the form
\begin{equation} \label{gs}
\frac{dx}{dt}
    =  - \sum_{(p, q) \in  S} a_{p q}x^{p+1}{y}^{q}
     =   P(x,y), \ \
\frac{dy}{dt}
    =  \sum_{(p,q) \in  S}
b_{qp}x^{q}{y}^{p+1}
     =  Q(x,y),
\end{equation}
where $ S$ is the set
$
 S = \{(p_j,q_j) \, | p_j+q_j \geq 0, j=1, \ldots, \ell \} \subset
(\{-1\} \cup \mathbb{N}_0 ) \times \mathbb{N}_0,
$
and $\mathbb{N}_0$ denotes the set of nonnegative integers. The
notation (\ref{gs}) simply emphasizes that we take into account only
nonzero coefficients of the polynomials and the symmetric form of the
system (in the sense that if there is a term $a_{pq}x^{p+1}y^q$ in the first equation of  (\ref{gs}), then there is also the term $b_{qp}x^{q}y^{p+1}$ in the second equation). 
Obviously, any polynomial system (\ref{eqgen}) can be
embedded into a system of the form (\ref{gs}).

We will assume that the parameters $a_{p_jq_j},\ b_{q_jp_j}$ ($j=1,
\ldots, \ell$) are from an infinite field $k$. Denote by $(a,b) = (a_{p_1 q_1},\dots,a_{p_\ell
q_\ell}, b_{q_\ell p_\ell} \dots, b_{q_1 p_1})$ the ordered vector
of coefficients of system (\ref{gs}), by $ E(a,b)$ the parameter
space of (\ref{gs}) (e.g. $ E(a,b)$  is $\mathbb{C}^{2\ell}$ or
$\mathbb{R}^{2\ell}$), and by $k[a,b]$ the polynomial ring in the
variables $a_{pq}$, $b_{qp}$ over the field $k$. It is clear that there is a
one--to--one correspondence between points of $E(a,b)$ and systems
of the form (\ref{gs}). The condition (\ref{crr1}) immediately
yields  that system (\ref{gs}) is time--reversible if and only if
\begin{equation} \label{ab} 
b_{qp} = \al^{p-q}a_{pq}, \qquad a_{pq} = b_{qp}\al^{q-p}.
\end{equation}
We rewrite (\ref{ab}) in the form
 \be \label{abp}
        a_{p_kq_k} = t_k, \quad
b_{q_kp_k} = \al^{p_k-q_k} t_k
\ee
  for $k =  1, \ldots, \ell$. From a geometrical point of view equations (\ref{abp}) define a surface
  in the affine space
  $
  k^{3
  \ell+1}=(a_{p_1q_1},\dots,a_{p_\ell q_\ell},b_{q_\ell p_\ell},\dots,b_{q_1p_1},t_1,\dots,t_\ell,\al). $
  Thus the set of all time-reversible systems is the projection of
  this surface onto  $k^{2\ell}=E(a,b)$. Nowadays  methods to compute such
  projections are well known (see e.g. \cite[Chapter 3]{Cox}). Namely, let 
  \be \label{H}
 H= \langle 1-w \,\tilde  \al_1\cdots\tilde \al_\ell, \
a_{p_kq_k} - t_k, \ \ \tilde \al_{k} b_{q_kp_k}-
\tilde{\tilde{ \al}}_{k} t_k \mid k=
1,\ldots ,\ell\rangle, \ee
where $\tilde{ \al}_{k}= \al^{q_k-p_q},\ \tilde{\tilde{ \al}}_{k}= 1$ 
if $p_k-q_k\le 0 $,  
$\tilde{ \al_{k}}=1,\ \tilde{\tilde{ \al}}_k= \al^{p_k-q_k}$ if $p_k-q_k>0$
 and $\langle h_1,\dots,h_m\rangle$ denotes the ideal generated by polynomials $h_1,\dots,h_m$.   
 Then as a direct corollary of Theorem  2 of \cite[Chapter 3, \S 3]{Cox} we
immediately find the
  minimal variety 
 which contains the set of all time-reversible
  systems in the family (\ref{gs}) (we remind  that the {\it variety} of 
ideal $I$, denoted by $\vv(I)$, is the set of 
 common zeros of all polynomials of $I$; 
the {\it Zariski closure} $\overline W$ of a set $W$ in the affine space $k^n$  is the minimal variety which contains $W$). 
Namely, from the theorem 
\be \label{t.Z.cl} \overline{\mathcal{R}}=\vv(\mathcal{I})\ \ {\rm where}\
\ \mathcal{I}= k[a,b] \cap H,
\ee
that is, the Zariski closure  of the set $\mathcal{R}$ of
all time-reversible systems is the variety of the ideal $\mathcal{I}$.
Moreover, according to Theorem 2 of  \cite[Chapter 3, \S 1]{Cox}  in order to find a generating set for
the ideal $\mathcal{I}$ it is sufficient to compute a Groebner basis 
 \footnote{Groebner basis is a
cornerstone of many algorithms of computational algebra 
(see e.g. \cite{Cox} for the definitions). Mathematica, Maple etc.
 have routines to
compute it.} for $H$ 
with respect to an elimination order with $\{w,\al, t_k\}>\{ a_{p_kq_k},\
b_{q_kp_k}\}$  and take from the output list  those 
polynomials, which depend only on $
        a_{p_kq_k},\
b_{q_kp_k}\ (k =  1, \ldots, \ell)$.

We now will obtain another description of the ideal $\mathcal{I}$, which will link the time-reversibility to invariants of a subgroup of $SL(2,\C)$. 
Denote by  $\mathcal{M}$  the set of all solutions $\nu =
(\nu_1,\nu_2,\dots,\nu_{2l})$ with non--negative components of the
equation
\be \label{lsib}
\zeta_1 \nu_1 + \zeta_2 \nu_2 + \cdots + \zeta_\ell \nu_\ell +
\zeta_{\ell+1} \nu_{\ell+1} + \cdots +\zeta_{2\ell} \nu_{2\ell} = 0,
\ee
where $\zeta_j=p_j-q_j$ for $j=1,\dots,\ell$, $\zeta_j=
q_{2\ell-j+1}-p_{2\ell -j+1}$ for $j=\ell+1,\dots,2\ell$, that is,
$$ \zeta=
(p_1 -q_1, p_2-q_2, \dots, p_\ell-q_\ell, q_\ell-p_\ell,\dots,
q_1-p_1)
$$
(we remind that
$(p_j,q_j)$ are from the set $S$ defining  system (\ref{gs})).
 Obviously, $\mathcal{M}$ is an  Abelian monoid
(semigroup).
For $\nu =(\nu_1,\dots,\nu_{2\ell})\in \mathcal{M}$ we
 denote by $[\nu]$
the  monomial
\be \label{nu}
a_{p_1q_1}^{\nu_1}a_{p_2q_2}^{\nu_2}\cdots a_{p_\ell
q_\ell}^{\nu_\ell} b_{q_\ell
p_\ell}^{\nu_{\ell+1}}b_{q_{\ell-1}p_{\ell-1}}^{\nu_{\ell+2}}\cdots
b_{q_1p_1}^{\nu_{2\ell}}
\ee and  by $\hat \nu$
{\it the involution of the vector} $\nu$,
$
\hat \nu= (\nu_{2\ell},\nu_{2\ell-1},\ldots ,\nu_1).
$
We call the ideal
\[
I_S = \langle [\nu]-[\hat \nu] \mid \nu\ \in \mathcal{M} \rangle
\subset k[a,b]
\]
 the \emph{Sibirsky ideal}  of system (\ref{gs})
(Sibirsky \cite{Sib1,Sib2} studied this ideal for the case of real systems
 embedded in the family (\ref{gs}); some studies in this direction were 
performed also in \cite{CGMM,Collins,Liu2}).
In the case that (\ref{gs}) is time-reversible, using (\ref{ab})
and (\ref{lsib}) we see that for $\nu \in \mathcal{M}$
\begin{equation}\label{cabr}
[\hat \nu] = \al^{\zeta \cdot \nu} [\nu] = [\nu],
\end{equation}
where $ \zeta \cdot \nu$ is the scalar product of $\zeta$ and $\nu$,
that is the left-hand side of (\ref{lsib}). By (\ref{cabr}) every
time--reversible system $(a,b) \in E(a,b)$ belongs to ${\bf V}(
I_S)$. However it is easily seen that the converse is false. The
following statement  proven in \cite{RS}  gives the exact
characterization of the set of time-reversible systems.
\begin{teo} \label{trev}
{
Let ${\mathcal{R}} \subset E(a,b)$ be the set of all
time--reversible systems in the family (\ref{gs}), then
\begin{enumerate}
\item[(a)] ${\mathcal{R}} \subset {\bf V}(I_S)$;
\item[(b)] ${\bf V}(I_S) \setminus {\mathcal{R}}
        = \{(a,b) \mid \exists (p,q) \in S\ \ {\rm such\ that}\
  \        a_{pq}b_{qp}=0\ \ {\rm  but }\ \ a_{pq}+b_{qp}\ne 0\}.$
\end{enumerate}}
\end{teo}
In other words, (b) means that if in a time-reversible  system
(\ref{gs}) $a_{pq}\ne 0$ then $b_{qp}\ne 0$ as well.
It follows from  (b) of the statement  that the inclusion in
 (a) is strict, that is ${\mathcal{R}} \subsetneqq {\bf
V}(I_S)$.

From  (\ref{t.Z.cl}) and (a) we have that $ \vv(\ic)\subset \vv(\is). $
We will   show  that in fact
 \be
\label{t.isym} I_S=\mathcal{I}\ {\rm and\ both\ ideals\ are\ prime}.
 \ee

When (\ref{t.isym}) is proven,  
from  (\ref{t.Z.cl}) and (\ref{t.isym}) we obtain 
the following description of time-reversible systems.
\begin{teo} \label{thm1} 
 {
The variety of the Sibirsky ideal $I_S $ is the Zariski closure of the set
$\mathcal{R}$ of all time-reversible systems in the family
(\ref{gs}).}
\end{teo}

We now 
  prove (\ref{t.isym}). 
 Suppose we are given the system 
\be\label{pse}
x_1 = \frac{f_1(t_1, \dots, t_m)}{g_1(t_1, \dots, t_m)}, \
      \dots, \
x_n = \frac{f_n(t_1, \dots, t_m)}{g_n(t_1, \dots, t_m)},
\ee where $f_j, g_j \in k[t_1, \dots, t_m]$ for $j = 1, \dots, n$.
Let $k(t_1, \dots, t_m)$ denote the ring of rational functions  in
$m$ variable with coefficients in $k$, and consider the ring
homomorphism
\[
\tilde \psi : k [\xn, t_1, \dots, t_m, w] \to k(t_1, \dots, t_m)
\]
defined by
$
t_i \to t_i, \quad x_j \to f_j(t_1, \dots, t_m) / g_j(t_1, \dots,
t_m),  w \to 1/g(t_1, \dots, t_m),
$
$i = 1, \dots, m$, $j = 1, \dots, n$ and $g=g_1g_2\cdots g_n$. Let
$$
\tilde H = \langle
    1 - w g,
    x_1 g_1(t_1, \dots, t_m) - f_1(t_1, \dots, t_m),
    \dots,\\
    x_n g_n(t_1, \dots, t_m) - f_n(t_1, \dots, t_m)
\rangle. 
$$
 It is not difficult to check that 
  \be \label{Hker}\tilde H = \ker(\tilde
\psi). \ee 
 Since $k[\xn, t_1, \dots, t_m, w] $ is a domain  (\ref{Hker}) yields that  $ \tilde H $ is a prime
ideal.

 By (\ref{Hker}) the ideal $ H$ defined by (\ref{H}) is the kernel of the ring homomorphism
$$
\psi: k[a,b,t_1,\dots,t_\ell,\al,w] \longrightarrow
k(\al,t_1,\dots,t_\ell)
$$
defined by $
        a_{p_kq_k} \mapsto t_k, \
b_{q_kp_k} \mapsto  \al^{p_k-q_k} t_k,\
w\mapsto 1/(\tilde \al_1\cdots\tilde \al_\ell)
$
for $k =  1, \ldots, \ell$.
We obtain a reduced Groebner basis $G$ of
$k[a,b]\cap H$ by computing a reduced Groebner basis of $H$ using an
elimination ordering with $\{a_{p_jq_j},\ b_{q_jp_j}\}<\{w, \al, t_j\} $ 
for all
$j=1,\dots, \ell$, and then intersecting it with $k[a,b]$.  Since
$H$ is binomial,  any reduced Groebner basis $G$ of $H$ also
consists of binomials. This shows that $H$ is a
binomial ideal.

It is easily  seen that any binomial of the form $[\alpha]-[\hat
\alpha] $ ($\alpha=(\alpha_1,\dots,\alpha_{2\ell})$ and $[\alpha]$ is defined by (\ref{nu})) is in $H$ (which is equal to $\ker(\psi)$) yielding $[\alpha]-[\hat \alpha]\in
k[a,b]\cap H=\mathcal{I} $ and, therefore, $I_S\subset \mathcal{I}.$
We  prove now that \be \label{iis} \mathcal{I}\subset I_S.\ee
$\mathcal{I}$ is a binomial ideal, thus to verify (\ref{iis}) it is
sufficient to show that
 all  binomials of $\mathcal{I}$ belong also to $I_S$.
 Since $\mathcal{I}$ is prime it is sufficient to consider only
 irreducible binomials. Indeed if $[\alpha]-[\beta]\in \ic$ is not
 irreducible, then $[\alpha]-[\beta]= ( [\alpha_1]-[\beta_1])
 [\theta]$, where $[\alpha_1]-[\beta_1]$ is irreducible binomial. It is clear
 from the definition of $\psi$ that $H$ contains no monomials, thus, since $H$ is prime,
 $[\alpha_1]-[\beta_1]\in H$.
Let now  $[\alpha]-[\beta]$ be an irreducible binomial in $I$, that
is, ${\rm supp}(\alpha)\cap {\rm supp}(\beta)= \emptyset$. Then
\begin{eqnarray*}
\psi([\alpha] - [\beta]) =
 t_1^{\alpha_1} \cdots
t_\ell^{\alpha_\ell}t_\ell^{\alpha_{\ell+1}} \al ^{(p_\ell-q_\ell)
\alpha_{\ell+1}} \cdots t_1^{\alpha_{2\ell}} \al ^{(p_1-q_1)
\alpha_{\ell+1}} - 
\\
 t_1^{\beta_1} \cdots
t_\ell^{\beta_\ell}t_\ell^{\beta_{\ell+1}} \al ^{(p_\ell-q_\ell)
\beta_{\ell+1}} \cdots t_1^{\beta_{2\ell}} \al ^{(p_1-q_1)
\beta_{\ell+1}} =\\
 t_1^{\alpha_1+\alpha_{2 \ell}} t_2^{\alpha_2+\alpha_{2\ell-1}}\cdots
  t_\ell^{\alpha_\ell+\alpha_{l+1}} \al^{\zeta_{\ell+1}\alpha_{\ell+1}+
  \zeta_{ \ell+2}\alpha_{\ell+2}\dots+
  \zeta_{2\ell}\alpha_{2\ell}}
-\\
t_1^{\beta_1+\beta_{2 \ell}} t_2^{\beta_2+\beta_{2\ell-1}}\dots
t_\ell^{\beta_\ell+\beta_{l+1}} \al^{\zeta_{\ell+1}\beta_{\ell+1}+
  \zeta_{ \ell+2}\beta_{\ell+2}\dots+ \beta_{2\ell}\beta_{2\ell}}.
\end{eqnarray*}
Thus, $\psi([\alpha] - [\beta]) =  0$ if and only if
 \be \label{fst}
\alpha_i+\alpha_{2\ell-i+1} =  \beta_i+\beta_{2\ell-i+1} \mbox{ for all } i=
1,2,\dots,\ell
 \ee
and \be \label{sec} \zeta_{\ell+1}\alpha_{\ell+1}+
  \zeta_{ \ell+2}\alpha_{\ell+2}\dots+
  \zeta_{2\ell}\alpha_{2\ell}=\ \zeta_{\ell+1}\beta_{\ell+1}+
  \zeta_{ \ell+2}\beta_{\ell+2}\dots+ \zeta_{2\ell}\beta_{2\ell}.
\ee
Since $\rm{supp}(\alpha)\cap\rm{supp}(\beta)= \emptyset$,  (\ref{fst}) yields that
 \be \label{bal} \beta_i= \alpha_{2\ell-i+1} \quad {\rm for}\ i=1,\dots,2\ell. \ee
 Hence
$
\beta= (\beta_1,\ldots ,\beta_{2\ell})= (\alpha_{2\ell},\ldots ,\alpha_1)=
\hat {\alpha}.
$
Noting that $\zeta_i=-\zeta_{2\ell +1-i} $  and using (\ref{bal}) we
write (\ref{sec}) as
$
\zeta_{1}\alpha_{1}+
  \zeta_{2}\alpha_{2}\dots+
  \zeta_{\ell}\alpha_{\ell}+\zeta_{\ell+1}\alpha_{\ell+1}+
  \zeta_{ \ell+2}\alpha_{\ell+2}\dots+
  \zeta_{2\ell}\alpha_{2\ell}=0.
$
That means $\alpha$ satisfies the equation (\ref{lsib}) and, hence,
$\alpha\in\cal{M}$. This completes the proof of (\ref{t.isym}).

Theorem 2 of \cite[Chapter 3, \S 1]{Cox} provides the following
algorithm for computing a generating set for the ideal $\mathcal{I}$ and, therefore, for the ideal ${I_S}$.
\begin{itemize}\item
Compute a Groebner basis  $G_H$ for $H$ defined by (\ref{H}) with
respect to any elimination order   with $\{w,\al,t_k\}
> \{{a_{p_kq_k}, b_{q_kp_k}} \mid  k=1,\dots,\ell\}$;
 \item the set $G_H\cap k[a,b] $ is a generating set for
 $\mathcal{I}$;
\item
 $ {\cal H}=  \{ \nu, \bar{\nu} \mid [\nu] - [\bar{\nu}] \in G\}
\cup  \{ {\bf e_i +e_j} \mid j =  2\ell-i+1; i = 1,\ldots, \ell\}$
is a Hilbert basis for the monoid $\cal M$ (here ${\bf e}_i$ is the
element $(0,\ldots ,0,1,0,\ldots ,0)$, with the non-zero entry in
the $i$-th position).
\end{itemize}
The correctness of the last step of the algorithm follows from
Theorem 4  of \cite{JRL}.

As an  example consider the  real system
\be \label{qvr}
\begin{array}{l} 
{\dot u} = -v+  a_{1} u^2+ a_{2} u v +a_{3}v^2 ,\quad 
{\dot v}= u+ b_{1}u^2 +b_{2} u v +b_{3}v^2.
\end{array}
\ee

Using the complexification procedure described in Section 2 
   we obtain from
(\ref{qvr})  the system
 \be\label{qv}
\begin{array}{l}
{\dot x} =i x +a_{10}x^2 +a_{01}xy +a_{-12}y^2 =ix+\widetilde P_2(x,y),\\
{\dot y}= -iy +b_{10}x y +b_{01}y^2 +b_{2,-1}x^2=-iy+\widetilde
Q_2(x,y).
\end{array}
\ee
 In the particular case $b_{ks}=\bar a_{ks}, \ y=\bar x$,
$a_{10}=
 ( ( a_1  +  b_2-a_3 )+i( b_1- a_2    - b_3))/4,$ 
$a_{01}= (a_1+a_3+i (b_1+b_3) )/2 $,
 $a_{-12}= (a_1-a_3  -  b_2 +i(   b_1+ a_2    - b_3))/4$
  system (\ref{qv}) is equivalent to
(\ref{qvr}).

 Computing a Groebner
basis of the ideal
$$
\mathcal{J}=\langle 1-w \al^4,
   {a_{10}} -
 {t_1},
   {b_{01}} - \al  {t_1},
   {a_{01}} -  t_2,\ 
  \al  {b_{10}} -  t_2,
   {a_{-12}} -  t_3,
   {\al^3b_{2,-1}} -  t_3    \rangle
$$
with respect to the  lexicographic order with $ w>\al>{t_1} > {t_2}
> {t_3}> a_{10}> a_{01}>a_{-12}>b_{10}> b_{01}>b_{2,-1} $ we obtain
a list of polynomials. (The list is too long to be presented here,
but is easily computed using any computer algebra system.) According
to step 2 of the algorithm  
in order 
 to obtain  a basis of $I_S$
we just have to pick up the polynomials that do not depend on $ w,
\al,{t_1} , {t_2} , {t_3}. $ 
In fact, 
there are five such
polynomials in the list:\
 $f_1= a_{01}^3 b_{2,-1}-a_{-12} b_{10}^3$, $ $
  $f_2= a_{10} a_{01} - b_{01} b_{10}$,\
  $f_3= a_{10}^3 a_{-12}- b_{2,-1} b_{01}^3$,\
  $f_4= a_{10} a_{-12} b_{10}^2- a_{01}^2 b_{2,-1} b_{01}$,\
  $f_5= a_{10}^2 a_{-12} b_{10} - a_{01} b_{2,-1} b_{01}^2.$
Thus, for system (\ref{qv})  
\[\is=\langle f_1,\dots,f_5\rangle.\]


An algorithm for computing generators of the 
 ideal $I_S$ has been also obtained  in
\cite{JRL}. To compare the efficiency of the algorithms  we carried out
 computations for few families of   system (\ref{gs}) on a PC with a 1.7 GHz
processor and   1.5 Gb RAM  using  Mathematica 4. 
It takes 0.2 seconds   CPU time  to find $I_S$ 
 using the algorithm described above  and 1.7 s CPU using the algorithm 
of \cite{JRL}. Replacing in (\ref{qv}) the polynomials $\widetilde P_2$ and $ \widetilde Q_2$ 
by the homogeneous   polynomials of the third
 and fourth degrees we find that it takes, respectively, 
0.04 and 4 s CPU with the algorithm of the present paper, and 
480 and 33450 s CPU with the algorithm of \cite{JRL}.

\section{ Time-reversibility and invariants}

We now show an interconnection of   time-reversibility and invariants
of a group of transformations of the phase space of system (\ref{gs}).
 Let $G$ be a matrix
group acting on ${\bf x} = (x_1, \dots, x_n)$ and let $k$ be any
field. We recall that a  polynomial $f({\bf x}) \in k[{\bf x}]$ is
\emph{invariant under $G$} if $f({\bf x}) = f(A \cdot {\bf x})$ for
every $A \in G$.
Consider  the  transformations  of the phase
space of (\ref{gs})
\begin{equation} \label{ROT}
  x'= \eta x, \quad  y'=\eta^{-1} y \ \ (x,y, \eta \in \C,\ \eta \ne 0).
\end{equation}
 The transformations (\ref{ROT}) form 
a subgroup of $SL(2,\C)$.   In $(x',y')$ coordinates system
(\ref{gs}) has the form
\[
\dot x' = \sum_{(p,q) \in   S} a(\eta )_{(p,q)}x'^{p+1}{y'}^{q}, \
\dot y' = \sum_{(p,q) \in  S} b(\eta)_{(q,p)}x'^{q}{y'}^{p+1}
\]
and the coefficients of the transformed system are
\begin{equation} \label{tr}
a(\eta)_{p_kq_k} = a_{p_kq_k} \eta^{q_k-p_k}, \quad b(\eta)_{q_kp_k}
= b_{q_kp_k} \eta^{p_k-q_k},
\end{equation}
where $k = 1, \dots, \ell$. Let $U_\eta$ denote the transformation
(\ref{tr}). 
$U_\eta$ is a
linear representation of  group (\ref{ROT}) in $\C^{2\ell}$.    We will write (\ref{tr}) in the short form
$
(a(\eta),b(\eta))  = U_\eta (a,b).
$
It is straightforward to see that a polynomial $f(a,b) \in
\mathbb{C}[a,b]$ is an invariant of the group ${{U_\eta}}$ if and
only if each of its terms is an invariant.

The action of $U_\eta$ on the coefficients $a_{ij}, b_{ji}$ of the
system of differential equations (\ref{gs}) yields the following
transformation of the monomial $[\nu]$ defined by (\ref{nu}):
\begin{eqnarray} \label{in59}
U_\eta[\nu] ={a(\eta)}_{p_1q_1}^{\nu_1} \cdots
                            {a(\eta)}_{p_\ell q_\ell}^{\nu_\ell} \
                            {b(\eta)}_{q_\ell p_\ell}^{\nu_{\ell+1}} \cdots
                            {b(\eta)}_{q_1p_1}^{\nu_{2\ell}}
                          = \\  \eta^{\zeta \cdot \nu}
                            a_{p_1q_1}^{\nu_1} \cdots
                            a_{p_\ell q_\ell}^{\nu_\ell} \
                            b_{q_\ell p_\ell}^{\nu_{\ell+1}} \cdots
                            b_{q_1p_1}^{\nu_{2\ell}}
                          =\eta^{\zeta\cdot \nu}[\nu]. \nonumber
\end{eqnarray}
Thus we see that \emph{ the monomial $[\nu]$ is invariant under the
action of  ${U}_\eta$ if and only if $\zeta\cdot \nu=0$, i.e.,
if and only if $\nu \in {{\mathcal{M}}}$.}

For $h \in \mathbb{C}^{2\ell}$   let $[\nu]|_{h}$
denote the value of the monomial $[\nu] = a_{p_1q_1}^{\nu_1} 
a_{p_2q_2}^{\nu_2} \cdots
b_{q_1p_1}^{\nu_{2\ell}}$ at the vector $h$, that is, $[\nu]|_{h} =
h_1^{\nu_1} h_2^{\nu_2}\cdots h_{2\ell}^{\nu_{2\ell}} $.
Denote by
 $\widehat{(a,b)}$  the involution of $(a,b)$, 
 \be \label{inv}
 \widehat{(a,b)} =
(b_{q_1p_1},\dots, b_{q_lp_l}, a_{p_lq_l},\dots, a_{p_1q_1}).
 \ee

Let $\mathcal{O}$  be an orbit of the group $U_\eta$ and assume that
there is $(a^*,b^*)\in \mathcal{O}$ such that
$\widehat{(a^*,b^*)}\in \mathcal{O}$. By
(\ref{ab}) and (\ref{tr})   the system
$(a^*,b^*)$ is time-reversible. We prove that   
then for any system $(a^0,b^0)\in
\mathcal{O}$ the system  $\widehat{(a^0,b^0)}$ is also in $\mathcal{O}$.
Let $(a^0,b^0)=U_{\eta_0}(a^*,b^*)$. Any invariant $[\nu]$
of the group $U_\eta$ is constant on any orbit of the group.
$(a^*,b^*)$ and $\widehat{(a^*,b^*)}$ belong to the same orbit $
\mathcal{O}$ of the group $U_\eta$, hence for all invariants $[\nu]$
of $U_\eta $ $[\nu]|_{(a^*,b^*)}= [\nu]|_{\widehat{(a^*,b^*)}}$.
Therefore, $[\nu]|_{(a^*,b^*)}= [\hat \nu]|_{(a^*,b^*)}$ yielding
$[\nu]|_{U_{\eta_0}(a^*,b^*)}= [\hat \nu]|_{U_{\eta_0}{(a^*,b^*)}}$
and, hence,
 \be \label{usl}
[\nu]|_{U_{\eta_0}(a^*,b^*)}= [ \nu]|_{\widehat
U_{\eta_0}{(a^*,b^*)}}, \ee where $\widehat U_{\eta_0}(a^*,b^*)$ is
the involution of the vector $ U_{\eta_0}(a^*,b^*)$. From
(\ref{usl}) and Theorem 3.2 of \cite{RS} we conclude that $\hat
U_{\eta_0}(a^*,b^*)$ belongs to the orbit $\mathcal{O}$ and,
therefore, the system $U_{\eta_0}(a^*,b^*)=(a^0,b^0)$ is time-reversible.

We say that the orbit $\mathcal{O}$ of the group $U_\eta$ is
invariant under the involution (\ref{inv}) if for any $(a,b)\in
\mathcal{O}$ the system $\widehat{(b,a)}$ also belongs to
$\mathcal{O}$.
By our reasoning above, 
 we have proven
\begin{teo}  (a)
The set of the orbits of $U_\eta$ is
divided into two not intersecting subsets: one consists of all
time-reversible systems and only time-reversible systems, and there
are no time-reversible systems in the other subset. 

(b)
The variety $\vv(\is)$ is the Zariski closure of all orbits of the
group $U_\eta$ invariant under the involution (\ref{inv}).
\end{teo}

To finish we note that the method described in the present paper 
can be applied to compute the set of systems 
 time-reversible with respect to some other 
symmetries, also in higher dimensions. However,  we believe that the
   beautiful
algebraic structure exhibiting in the parameter space  by the symmetry (\ref{xy}) is due to the fact that (\ref{xy}) is an involution.  
Note also that an interesting future project would be a generalization of the 
presented results to the case of higher dimensional systems.

\section*{Acknowledgements}

This   work was supported   by the Ministry of
Higher  Education, Science and Technology 
of the Republic of Slovenia,  Nova Kreditna Banka Maribor and 
 TELEKOM  Slovenije.

\small


\begin{thebibliography}{99}

\bibitem{Cairo} 
Cairo L,  Chavarriga J,  Gin\'e J, Llibre J 1999
{\em Comput. Math. Appl.} {\bf   38}  39--53

\bibitem{Chavarriga-1}
Chavarriga J, Gin\'e J and  Garc\'ia I A 2000
{\sl J. Computational and Applied Mathematics} 
{\bf 126} 351-368


\bibitem{Chavarriga-2}
Chavarriga J and Sabatini M
1999
{\sl Qual. Theory of Dyn. Systems} {\bf 1} 1--70


\bibitem{CGMM} Cima  A,  Gasull A,  Ma\~{n}osa V  and Manosas F 1997
               \emph{Rocky Mountain J.~Math.} {\bf 27}  471--501.

\bibitem{Cox}
   Cox D,   Little J and    O'Shea D 1992 
  {\it Ideals, Varieties,
   and Algorithms} (New York:Springer--Verlag)


\bibitem{Collins}  
Collins C~B 2001 
   {\it J.~Math.~Anal.~Appl.} {\bf 259}  168--187.

\bibitem{JRL}
  Jarrah A,    Laubenbacher R and  Romanovski V G
  2003 {\it
    J.\ Symb.\ Comput.} {\bf
 35} 577--589


\bibitem{Lamb}
   Lamb J S W and   Roberts J A G
  1998 {\it
  Phys.\ D}
{\bf
  112}
  1--39
\bibitem{Liu2} Li Y--R, and  Li J--B 1989 Theory of values of singular point 
   in complex autonomous differential systems. {\it Sci.~China Series A} 
   {\bf 33}  10--23

\bibitem{MMV}
Marde\v si\'c P,  Marin D,  Villadelprat J 2006 {\em  J. Differential 
Equations}  {\bf 224}   120--171 

\bibitem{Valery-1}
Romanovski V G 2003 
{\sl Progr. Theoret. Phys. Suppl.} {\bf 150} 243--254

\bibitem{RS}
  Romanovski V G and  Shafer D S 2005
{\it Differential Equations with Symbolic Computations. Trends in Mathematics}  (Eds  Wang D and  Zheng Z,
  Birkhauser)
  67-84

\bibitem{Sib1} Sibirsky K S 1976 \emph{  Algebraic Invariants of Differential 
               Equations and Matrices.}
                Kishinev: Shtiintsa (in Russian)
\bibitem{Sib2} 
 Sibirsky K S 1982 {\it
Introduction to the Algebraic Theory of
   Invariants of  Differential Equations}
(Kishinev:
  Shtiintsa) (in Russian,
   English transl.:  Manchester University Press, 1988)



\end{thebibliography}
\end{document}